\documentclass[12pt]{article}
\usepackage[cp1251]{inputenc}

\usepackage[english]{babel}
\usepackage{amsmath,amsfonts,amscd,amssymb,latexsym}

\makeindex
\newtheorem{theorem}{Theorem}
\newtheorem{proposition}{Proposition}

\newtheorem{corollary}{Corollary}
\newtheorem{definition}{Definition}
\newtheorem{example}{Example}
\newtheorem{remark}{Remark}

\begin{document}


\begin{center} {\Large Rota-Baxter operators and non-skew-symmetric solutions of the classical Yang-Baxter equation on quadratic Lie algebras.}
\end{center}

\vspace{5mm}

\begin{center}
{\bf M.\,E.\,Goncharov}
\end{center}

\begin{center}
Sobolev Institute of Mathematics\\
goncharov.gme@gmail.com
\end{center}

\begin{abstract}
We study possible connections between Rota-Baxter operators of
non-zero weight and non-skew-symmetric solutions of the classical
Yang-Baxter equation on finite-dimensional quadratic Lie algebras.
The particular attention is made to the case when for a solution $r$
the element $r+\tau(r)$ is $L$-invariant.
\end{abstract}

\section{Introduction}

Let $A$ be an arbitrary algebra over a field $F$, $\lambda\in F$. A
map $R:A\mapsto A$ is called a Rota-Baxter operator of weight
$\lambda$ if for all $x,y\in A$:
\begin{equation}\label{e1}
R(x)R(y)=R(R(x)y+xR(y)+\lambda xy)
\end{equation}

As an example of a Rota-Baxter operator of weight zero one can
consider the operation of integration on the algebra of continuous
functions on $\mathbb{R}$: the equation \eqref{e1} follows from the
integration by parts formula.

Rota-Baxter operators for associative algebras first appears in the
paper of G. Baxter as a tool for studying integral operators that
appears in the theory of probability and mathematical statistics
\cite{Br}. The combinatorial properties of Rota---Baxter algebras
and operators were studied in papers of F.V. Atkinson, P. Cartier,
G.-C. Rota and the others (see \cite{Atk}-\cite{Car}). For basic
results and the main properties of Rota---Baxter algebras see
\cite{Guo}.

Independently, in early 80-th Rota-Baxter  operators on Lie algebras
naturally appears in papers of A.A. Belavin, V.G. Drinfeld \cite{BD}
and M.A. Semenov-Tyan-Shanskii \cite{STS} while studying the
solutions of the classical Yang-Baxter equation.

There is a standard method for constructing Rota---Baxter operations
of weight zero on a quadratic (that is, possessing a non-degenerate
 invariant form $\omega$) Lie algebra $(L,\omega)$ from
skew-symmetric solutions of the classical Yang---Baxter equations
(CYBE): if $r=\sum a_i\otimes b_i$ is a skew-symmetric solution of
CYBE (that is, $\tau(r)=-r$ where $\tau$ is the switch morphism),
then one can define an operator $R$ on $L$ by
$$R(a)=\sum\omega(
b_i,a) a_i.$$
 It turns out that $R$ is a Rota---Baxter operator of weight
 zero\cite{BD,STS}.
  Moreover, Rota-Baxter operators of weight zero that can be obtained from skew-symmetric solutions of CYBE can be easily described:
   they satisfy $R+R^*=0$ where $R^*$ is the adjoint to $R$ with respect to the form $\omega$ operator.

 The case when $r$ is a non-skew-symmetric solution of
the classical Yang-Baxter equation was considered in \cite{GME}. It
was proved that if $L$ is a simple Lie algebra and $r$ is a solution
of CYBE such that the element $r+\tau(r)$ is $L$-invariant, then
there is a non-degenerate invariant form $\omega$ on $L$ such that
the corresponding to $r$ linear map $R$ is  a Rota-Baxter operator
of a non-zero weight $\lambda$. In this paper  we complement this
result by proving that the obtained operator $R$ satisfies
$R+R^*+\lambda id=0$, where $\lambda$ is the weight of $R$ and $id$
is the identity operator.

There is a connection between solutions of CYBE on Lie algebras and
Lie bialgebras that were introduced by Drinfeld \cite{Drinf} for
studying the solutions of the classical Yang---Baxter equation on
Lie algebras.  If $L$ is a Lie algebra and  $r=\sum\limits_i
a_i\otimes b_i\in L\otimes L$, then one can define a
comultiplication $\delta_r: L\mapsto L\otimes L$ as
$$
\delta_r(a)=\sum [a_i,a]\otimes b_i+a_i\otimes [b_i,b].
$$
If the element $r$ is a solution of CYBE and $r+\tau(r)$ is
$L$-invariant, then the pair $(L,\delta_r)$ is a Lie bialgebra.
Bialgebras of this type are called quasitriangular (or triangular,
if $r$ is skew-symmetric). Triangular Lie bialgebras play an
important role as they lead to solutions of quantum Yang-Baxter
equation.

In the current paper, we consider the correspondence between
solutions of CYBE and Rota-Baxter operators of non-zero weight on
arbitrary quadratic  Lie algebra $(L,\omega)$. In section 3 we
describe Rota-Baxter operators that can be obtained from arbitrary
solutions of CYBE on arbitrary  (not necessarily simple) quadratic
finite-dimensional Lie algebra. Unlike the case of simple Lie
algebras, here $R+R^*+\lambda id$ is not necessarily equal to zero.

In section 4 we study the connection between Rota-Baxter operators
of nonzero weights and solutions $r$ of CYBE for which $r+\tau(r)$
is $L$-invariant.  It turns out that the corresponding Rota-Baxter
operators  also do not necessarily satisfy the equality
$R+R^*+\lambda id=0$. Nevertheless, in this case for such a
Rota-Baxter operator $R$ there is an ideal $I\subset [L,L]$ such
that $R(I)=0$ and the restriction of $R$ on the quotient algebra
$L/I$ satisfy $R([a,b])+R^*([a,b])+\lambda [a,b]=0$ for all $a,b\in
L/I$.

Also in section 4 we consider the situation when for a solution $r$
of CYBE the map $R$ and the adjoint map $R^*$ are Rota-Baxter
operator of the same nonzero weight.

\section{Definitions and preliminary results}

Throughout the paper the characteristic of the ground field $F$ is 0
and all spaces are supposed to be finite-dimensional.

 Given vector spaces $V$ and $U$ over a field $F$, denote
by $V\otimes U$ its tensor product over $F$. Define the linear
mapping $\tau$ on $V$ by $\tau(\sum\limits_ia_i\otimes
b_i)=\sum\limits_ib_i\otimes a_i$. Denote by $V^*$ the dual space of
$V$. For $f\in V^*$ and $v\in V$ by $\langle f,v\rangle$ we will
denote the action of $f$ on the vector $v$, that is $\langle
f,v\rangle=f(v)$.

Let $L$ be a Lie algebra with a multiplication $[\cdot,\cdot]$. By
$Z(L)$ we will denote the center of $L$. Recall, that $L$ acts on
$L^{\otimes n}$ by
$$
[x_1\otimes x_2\otimes \ldots \otimes x_n,y]=\sum\limits_i
x_1\otimes\ldots\otimes [x_i,y]\otimes \ldots\otimes x_n$$

for all $x_i,y\in L$.

\begin{definition} An element $r\in L^{\otimes n}$ is called
$L$-invariant (or $ad$-invariant) if $[r,y]=0$ for all $y\in L$.
\end{definition}

For the special linear Lie algebra $sl_2(F)$ we will denote by
$e,h,f$ the standard basis of $sl_2(F)$: $e=e_{12}$,
$h=e_{11}-e_{22}$, $f=e_{21}$.

The next definition of Lie bialgebra that was given in \cite{Drinf}.

\begin{definition} Let $L$ be a Lie algebra with a comultiplication
$\delta$. The pair $(L,\delta)$ is called a Lie bialgebra if and
only if $(L,\delta)$ is a Lie coalgebra and $\delta$ is a 1-cocycle,
i.e., it satisfies
$$
\Delta([a,b])=\sum([a_{(1)},b]\otimes
a_{(2)}+a_{(1)}\otimes[a_{(2)},b])+\sum([a,b_{(1)}]\otimes
b_{(2)}+b_{(1)}\otimes[a,b_{(2)}])
$$
for all $a,b\in L$.
\end{definition}

There is an important type of Lie bialgebras called coboundary Lie
bialgebras. Namely, let $L$ be a Lie algebra and
$r=\sum\limits_ia_i\otimes b_i\in L\otimes L$. Define a
comultiplication $\delta_r$ on $L$ by
$$
\Delta_r(a)=[r,a]=\sum\limits_i[a_i,a]\otimes b_i+a_i\otimes[b_i,a]
$$
for all $a\in L$. It is easy to see that $\delta_r$ is a 1-cocycle.

Define an element $C_L(r)$ as
$$
C_L(r)=[r_{12},r_{13}]-[r_{23},r_{12}]+[r_{13},r_{23}].
$$
Here $[r_{12},r_{13}]=\sum\limits_{ij}[a_i,a_j]\otimes b_i\otimes
b_j$, $[r_{23},r_{12}]=\sum\limits_{ij}a_i\otimes[a_j,b_i]\otimes
b_j$,
 and $[r_{13},r_{23}]=\sum\limits_{ij} a_i\otimes a_j\otimes [b_i,b_j]$.

 The dual algebra $L^*$ of the coalgebra
$(L,\delta_r)$ is anti-symmetric if and only if $r+\tau(r)$ is
$L$-invariant. Also, $L$  satisfies the Jacobi identity if and only
if $C_L(r)$ is $L$-invariant. In particular, if $r=\sum a_i\otimes
b_i\in L\otimes L$ satisfy
\begin{equation}\label{lieYB}
\sum\limits_{ij}[a_i,a_j]\otimes b_i\otimes
b_j-a_i\otimes[a_j,b_i]\otimes b_j+a_i\otimes a_j\otimes
[b_i,b_j]=0,
\end{equation}
and $r+\tau(r)$ is $L$-invariant, then $(L,\delta_r)$ is a Lie
bialgebra called quasitriangular Lie bialgebra (or triangular, if
$\tau(r)=-r$). The equation \eqref{lieYB} is called \emph{the
classical Yang---Baxter
 equation}.

Note that the equation \eqref{lieYB},  as well as the corresponding
to it bialgebra structures, can be considered for every variety of
algebras. For Jordan, associative, alternative and Malcev bialgebras
it is known that if $r$ is a skew-symmetric solution of the
classical Yang---Baxter equation on an algebra $A$, then
$(A,\delta_r)$ is a bialgebra of corresponding variety
\cite{Zhelyabin98,Zhelyabin,versh,Aquiar,Polishchuk,Gme,GMM}. Also,
 it is also worth mentioning  papers of C. Bai, L. Guo and X.
Ni\cite{BGN,BGN1} where it was studied the connection between
coboundary Lie bialgebras with a generalization of CYBE (extended
CYBE) and with a generalization of Rota-Baxter operators of nonzero
weight (extended $\mathcal O$-operators).

\begin{definition} Let $L$ be a Lie algebra. A bilinear form
$\omega$ on $L$ is called invariant if for all $a,b,c\in L$:
$\omega([a,b],c)=\omega(a,[b,c])$.
\end{definition}

\begin{definition} Let $L$ be a Lie algebra and $\omega$ is an
invariant non-degenerate form on $L$. Then the pair $(L,\omega)$ is
called a quadratic Lie algebra.
\end{definition}

Let $(L,\omega)$ be a quadratic finite-dimensional Lie algebra. Then
the dual space $L^*$ is naturally isomorphic to $L$ and we may
present $L^*$ as $L^*=\{ a^*|\ a\in L\}$, where $\langle
a^*,b\rangle=\omega(a,b)$ for all $a,b\in L$. Note that the
associativity of the form $\omega$ implies $\langle
[a,b]^*,c\rangle=\langle a^*,[b,c]\rangle=\langle b^*,[c,a]\rangle$.

Also, if $(L,\omega)$ is a quadratic finite-dimensional Lie algebra,
then there is a natural isomorphism between the space of
endomorphisms $End_F(L)$ and $L\otimes L$: for every $\varphi\in
End_F(L)$ there is a unique $r=\sum\limits_i a_i\otimes b_i\in
L\otimes L$ such that for all $a\in L$: $\varphi(a)=\sum\limits_i
\omega(b_i,a)a_i$.

There is a connection between skew-symmetric solutions of the
classical Yang-Baxter equation and Rota-Baxter operators on
quadratic Lie algebras: if $(L,\omega)$ is a quadratic Lie algebra
and $r=\sum\limits_i a_i\otimes b_i\in L\otimes L$ is a
skew-symmetric solution of the classical Yang-Baxter equation on
$L$, then an operator $R$ defined as
\begin{equation}\label{op1}
R(a)=\sum\limits_i \omega(b_i,a)a_i
\end{equation}
for all $a\in L$ is a Rota-Baxter operator of weight zero
\cite{BD,STS}. Moreover, if $R^*$ is the adjoint to $R$ with respect
to the form $\omega$ operator, then $R+R^*=0$.

In \cite{GME} it was proved that if $L$ is a simple  Lie algebra,
$r=\sum a_i\otimes b_i$ is a  non-skew-symmetric solution of CYBE
such that $\tau(r)+r$ is $L$-invariant, then there is a
non-degenerate bilinear invariant form $\omega$ on $L$ such that the
operator $R$ defined as in \eqref{op1} is a Rota-Baxter operator of
a non-zero weight.

If $L$ is not a simple quadratic Lie algebra and $r$ is a
non-skew-symmetric solution of CYBE then everything is possible: as
the following examples show, the corresponding map $R$ may be a
Rota-Baxter operator of zero or non-zero weight, or not a
Rota-Baxter operator at all.

\begin{example} Consider $L=gl_2(\mathbb{C})$ --- Lie algebra of
$2\times 2$ matrices over $F$ and $r=\frac{1}{2}(E\otimes
e_{11}+e_{22}\otimes E)$, where $E$ is the identity matrix. It is
easy to see that $r$ is a not skew-symmetric solution of CYBE and
$r+\tau(r)=E\otimes E$ is $L$-invariant. But the operator $R$
defined as is \eqref{op1} is a Rota-Baxter operator of weight zero.
\end{example}

\begin{example}  Let $r_i\in L_i\otimes L_i$ ($i=1,2$) be  a solution of CYBE on
simple complex Lie algebra $L_i$ such that $r_i+\tau(r_i)$ is
$L_i$-invariant. Let $L=L_1\oplus L_2$ be a semisimple complex Lie
algebra and define $r=r_1+r_2$. In this case $r$ is a solution of
CYBE on $L$, $r+\tau(r)$ is $L$-invariant and if at least one
$i\in\{1,2\}$ satisfies $r_i+\tau(r_i)\neq0$, then $r+\tau(r)\neq
0$.

It is known that every non-degenerate invariant form $\omega$ on $L$
may be presented as $\omega=\alpha_1 \chi_1+\alpha_2\chi_2$, where
$\chi_i$ is the Killing form on $L_i$ and $\alpha_i\neq 0$.

 In \cite{GME} it was proved that if $r_i+\tau(r_i)\neq 0$, then the operator $R_i$, defined as in
\eqref{op1} (where $\omega=\chi_i$) is a Rota-Baxter operator of
some nonzero weight $\lambda_i$. And if for $i=1,2$
$r_i+\tau(r_i)\neq 0$, then we can find nonzero scalars $\mu_i\in
\mathbb{C}$ ($i=1,2$) such that for the form $\omega=\mu_1
\chi_1+\mu_2\chi_2$ the operator $R$ defined by \eqref{op1} is a
Rota-Baxter operator of some nonzero weight $\lambda$.

But if $r_1+\tau(r_1)=0$ and $r_2+\tau(r_2)\neq 0$, then for every
form $\omega$ the restriction of the corresponding operator $R$ to
$L_1$ is a Rota-Baxter operator of weight zero while the restriction
of $R$ on $L_2$ will be a Rota-Baxter operator of some nonzero
weight $\lambda$. Therefore, in this case $R$ is not a Rota-Baxter
operator.

\end{example}

\section{CYBE and Rota-Baxter operators on quadratic Lie
algebras.}

Let $(L,\omega)$ be a quadratic Lie algebra and $r=\sum\limits_i
a_i\otimes b_i\in L\otimes L$.  Define an operator $R:L\mapsto L $
as in \eqref{op1}. By $R^*$ we denote the adjoint operator with
respect to the form $\omega$.

\begin{proposition} Let $r$ be a solution of the CYBE on $L$. Then the following equalities hold for all $x,y\in
L$:

1. $[R(x),R(y)]-R([x,R(y)])+R([R^*(x),y])=0$,

2. $[R^*(x),R^*(y)]+R^*([x,R(y)])-R^*([R^*(x),y])=0$,

3. $r+\tau(r)$ is L-invariant if and only if $R+R^*$ lies in  the
centroid of $L$, that is $$[(R+R^*)(x),y]=(R+R^*)([x,y]).$$
\end{proposition}
{\bf Proof}

1. For $x,y\in L$ consider a map $\psi_{x,y}:L\otimes L\otimes
L\mapsto L$ defined as
$$
\psi_{x,y}(a\otimes b\otimes c)=\omega(x,b)\omega(y,c)a
$$
for all $a,b,c\in L$.

 Since $r$ is a solution of CYBE, then for
all $x,y\in L$ we have $\psi_{x,y}(C_L(r))=0$. Direct computations
show
$$
0=\psi_{x,y}(C_L(r))=\sum\limits_{i,j}[a_i,a_j]\omega(b_i,x)\omega(b_j,y)-a_i\omega([a_j,b_i],x)\omega(b_j,y)+a_i\omega(a_j,x)\omega([b_i,b_j],y)=
$$
$$
= [R(x),R(y)]-R([x,R(y)])+R([R^*(x),y]).
$$

2. The proof is similar to 1.

3. Let $a,b\in L$ and consider a map $\psi_x=L\otimes L\mapsto L$
defined as
$$
\psi_x(a\otimes b)=\omega(x,a)b.
$$
Since the form $\omega$ is non-degenerate, an element $h\in L\otimes
L$ is equal to zero if and only if $\psi_x(h)=0$ for all $x\in L$.

 Now let $x,y\in L$.  Direct computation shows:

$$\psi_x([r+\tau(r),y])=\omega([a_i,y],x)b_i+\omega(a_i,x)[b_i,y]+\omega([b_i,y],x)a_i+\omega(b_i,x)[a_i,y]=$$
$$
=R^*([y,x])+[R^*(x),y]+R([y,x])+[R(x),y].
$$
And the proposition is proved.

The next theorem gives necessary and  sufficient conditions when $R$
is a Rota-Baxter operator of weight $\lambda$.

\begin{theorem} Let $(L,\omega)$ be a quadratic Lie algebra over a
field $F$ and $r=\sum a_i\otimes b_i\in L\otimes L$. Let $R$ be the
operator defined as in \eqref{op1} and $R^*$ be the adjoint to $R$
with respect to the form $\omega$ operator. Then

1. If $r$ is a solution of CYBE on $L$, then $R$ is a Rota-Baxter
operator of a weight $\lambda$ if and only if for all $a,b\in L$:

\begin{equation}\label{YBE-RB}
[R(a),b]+[R^*(a),b]+\lambda[a,b]\in ker(R).
\end{equation}
This is equivalent to:
$$
\omega(b_j,a)\omega([b_i,a_j],b)a_i+\omega(a_j,a)\omega([b_i,b_j],b)a_i+\lambda\omega(b_i,[a,b])a_i=0.
$$

2. Conversely, if $R$ is a Rota-Baxter operator of a weight
$\lambda$ on $L$, then  $r$ is a solution of CYBE on $L$  if and
only if the equality \eqref{YBE-RB} holds for all $a,b\in L$.
\end{theorem}

{\bf Proof} The first statement follows from the proposition 1.1. In
order to prove the second statement consider  a map
$\psi_{x,y,z}:L\otimes L\otimes L\mapsto L$ defined as
$\psi_{x,y,z}(a\otimes b\otimes
c)=\omega(x,a)\omega(y,b)\omega(z,c)$. It is well known that an
element $h\in L\otimes L\otimes L$ is equal to zero if and only if
for all $x,y,z\in L$: $\psi_{x,y,z}(h)=0$. Using the technique from
the proof of the proposition 1, we have
$$
\psi_{x,y,z}(C_L(r))=\omega([R(x),R(y)]-R([x,R(y)])+R([R^*(x),y]),z)=
$$
$$=\omega([R(x),R(y)]-R([x,R(y)])-R([R(x),y])-R(\lambda[x,y]),z)+
$$
$$+\omega(R([R^*(x),y]+R([R(x),y])+R(\lambda[x,y])),z)=
$$
$$=\omega(R([R^*(x),y]+R([R(x),y])+R(\lambda[x,y])),z).
$$

Thus, $C_L(r)=0$ if and only if for all $x,y\in L$:
$R([R^*(x),y])+R([R(x),y])+R(\lambda[x,y])=0$.

\begin{remark} The statement of the theorem 1 can be formulated for
$R^*$. The condition \eqref{YBE-RB} in this case should be replaced
by $$[R(a),b]+[R^*(a),b]+\lambda[a,b]\in ker(R^*).$$
\end{remark}

\begin{corollary}
Let $(L,\omega)$ be a quadratic Lie algebra and $R$ be a Rota-Baxter
operator of a nonzero weight $\lambda$. If for all $a\in L$:
$$
R(a)+R^*(a)+\lambda a\in Z(L),
$$
then the corresponding tensor $r$ is a solution of CYBE.
\end{corollary}

\begin{example} Let $L=sl_2(F)$ and $\chi$ be the Killing form on
$L$, $F[t]$ be the algebra of polynomials in the variable $t$ and
$F_2[t]=F[t]/(t^2)$ be the quotient algebra over the ideal spanned
by $t^2$. Define $L_2$ as the tensor product $L\otimes_F F_2[t]$.
The product on $L_2$ is defined as
$$
[a\otimes \overline{f(t)},b\otimes \overline{g(t)}]=[a,b]\otimes
\overline{f(t)g(t)}
$$
where $[a,b]$ is the product in $L$ and $\bar{\ }:F[t]\mapsto
F[t]/(t^2)$ is the natural homomorphism. On $L_2$ consider a form
$\omega$ defined as
$$
\omega(a\otimes \overline{f(t)},b\otimes
\overline{g(t))}=\chi(a,b)\pi(f(t)g(t)),
$$
where $\pi(f(t))=f_0+f_1$ whenever $f(t)=f_0+f_1t+\ldots+f_nt^n$.

Then $(L_2,\omega)$ is quadratic Lie algebra \cite{BajoBen}. We can
identify elements $a\otimes \overline{1}$ with $a$ and $a\otimes
\overline{t}$ with $a\overline{t}$ and present $L_2$ as sum of the
subalgebras $L$ and $L\overline{t}$: $L_2=L\oplus L\overline{t}$.
Define an operator $R:L_2\mapsto L_2$ as
$$
R(a)=a,\ R(a\overline{t})=0
$$
for all $a\in L$.

It is known that $R$ is a Rota-Baxter operator of weight $-1$ (see
\cite{Guo}). Direct computations show that for all $a\in L$:
$R^*(a)=R^*(a\overline{t})=a\overline{t}$. We have:
$$
R(a)+R^*(a)-a=a\overline{t},\ \
R(a\overline{t})+R^*(a\overline{t})-a\overline{t}=0.
$$
It means that for all $x\in L_2$:
$[R(x)+R^*(x)-x,L_2]=L\overline{t}\subseteq ker(R)$. The map $R$ is
defined by the tensor $r=\frac{1}{4}h\otimes h\overline{t}+e\otimes
f\overline{t}+f\otimes e\overline{t}$. By theorem 1 $r$ is the
solution of CYBE.

Note, that $R^*$ is not a Rota-Baxter operator. And in order to
prove this now it is not necessary to check \eqref{e1}. Indeed, we
know, that $r$ is a solution of CYBE and by theorem 1 $R^*$ is a
Rota-Baxter operator of a weight $\lambda$ if and only if  for all
$x,y\in L_2$ $[R(x),y]+[R^*(x),y]+\lambda[x,y]\in ker(R^*)$. We have
already proved that $[R(x)+R^*(x)-x,L_2]=L\overline{t}$ so $R^*$ is
not a Rota-Baxter operator of weight $-1$. Let $\lambda\neq -1$.
Then for $a,b\in L$:
$$
R^*([R(a\overline{t})+R^*(a\overline{t})+\lambda
a\overline{t},b])=(1+\lambda)R^*([a,b]\overline{t})=(1+\lambda)[a,b]\overline{t}
$$
that proves that $R^*$ is not a Rota-Baxter operator.
\end{example}

\begin{example} Let $(L_2,\omega)$ be the quadratic Lie algebra
from example 3. Define a map $Q$ as the projection on the subalgebra
$L\overline{t}$ with the kernel $ker(Q)=L$. Then $Q$ is a
Rota-Baxter operator of weight $-1$. Note, that $Q=id - R$, where
$id$ is the identity operator and $R$ is the Rota-Baxter operator
from the example 3. It means that $Q^*=id-R^*$, that is for all
$a,b\in L$: $Q^*(a)=a-a\overline{t}$, $Q^*(a\overline{t})=0$. But
for all $a,b\in L$:
$$
Q([Q(a)+Q^*(a)-a,b])=Q([a,b]-[a,b]\overline{t})=-[a,b]\overline{t}.
$$
And by theorem 1 an element  $h=\frac{1}{4}h\overline{t}\otimes
(h-h\overline{t})+x\overline{t}\otimes
(y-y\overline{t})+y\overline{t}\otimes(x- x\overline{t})$ is not a
solution of CYBE on $L_2$.
\end{example}

\section{Solutions of CYBE with $L$-invariant symmetric part and Rota-Baxter operators on quadratic Lie
algebras.}

In this section, we will consider possible connections between
solutions of CYBE such that $r+\tau(r)$ is  $L$-invariant and
Rota-Baxter operators on quadratic Lie algebra $(L,\omega)$. Note,
that it is known that skew-symmetric solutions of CYBE are in one to
one correspondence with Rota-Baxter operators $R$ of weight zero
satisfying $R+R^*=0$ (see \cite{BD,STS}). We will study the case
when $r+\tau(r)\neq 0$.

Unless otherwise is specified, throughout this section $(L,\omega)$
is a finite-dimensional quadratic Lie algebra over a field $F$,
$r=\sum a_i\otimes b_i\in L\otimes L$, $R$ is the operator defined
as in \eqref{op1} and $R^*$ is the adjoint to $R$ with respect to
the form $\omega$ operator.

{\bf Definition} For a linear map $R:L\mapsto L$ and a scalar
$\alpha\in F$ define an operator
$$
\theta_{\alpha}=R+R^*+\alpha id,
$$
where $id$ is the identity operator.

\begin{theorem} Let $r=\sum a_i\otimes b_i\in L\otimes L$ be a
solution of CYBE on $L$ such that $r+\tau(r)$ is $L$-invariant.
Then, for every $\lambda$, a set
\begin{equation}\label{ideal}
I_{\lambda}=\{\theta_{\lambda}(x)=R(x)+R^*(x)+\lambda x|\ \ x\in
[L,L]\}
\end{equation}
 is an ideal in $L$. If
$I_\lambda$ is $R$-invariant, then $I_{\lambda}$ is also
$R^*$-invariant and the restrictions of $R$  and $R^*$ on the
quotient algebra $L/I_{\lambda}$ are Rota-Baxter operators of weight
$\lambda$. Moreover, in this case we have :
\begin{equation}\label{con1}
\theta_{\lambda}(x)=R(x)+R^*(x)+\lambda x=0\ \text{for all  $x\in
[L/I_{\lambda},L/I_{\lambda}]$.}
\end{equation}
\end{theorem}

{\bf Proof} Consider $\lambda\in F$ and define
$$
I_{\lambda}=\{R(x)+R^*(x)+\lambda x|\ \ x\in [L,L]\}.
$$
 For all $z\in L$, using proposition
1 (3), we have:
$$
[R([x,y])+R^*([x,y])+\lambda
[x,y],z]=[R([x,y]),z]+[R^*([x,y]),z]+\lambda [[x,y],z]= $$
$$=R([[x,y],z])+R^*([[x,y],z])+\lambda [[x,y],z]\in I_{\lambda}.$$

Thus, $I_\lambda$ is an ideal in $L$. Note, that by proposition
1(3):
$$
R[x,y]+R^*([x,y])+\lambda [x,y]=[R(x)+R^*(x)+\lambda x,y]\subset
[L,L].
$$

 Suppose that
$R(I_{\lambda})\subset I_{\lambda}$. First we prove that
$I_{\lambda}$ is also an $R^*$-invariant ideal. Consider
$\theta_{\lambda}([x,y])\in I_{\lambda}$. We have

$$
R(\theta_{\lambda}([x,y]))+R^*(\theta_{\lambda}([x,y]))=R([R(x)+R^*(x)+\lambda
x,y])+R^*([R(x)+R^*(x)+\lambda x,y])=
$$
$$
=R([R(x)+R^*(x),y])+R^*([R(x)+R^*(x),y])+\lambda
(R([x,y])+R^*([x,y]))=
$$
$$
=R([R(x)+R^*(x),y])+R^*([R(x)+R^*(x),y])+\lambda [R(x)+R^*(x),y]\in
I_{\lambda}.
$$
 And since $R(\theta_{\lambda}([x,y]))\in I_{\lambda}$, we have that $R^*(\theta_{\lambda}([x,y]))\in
 I_{\lambda}$.

Consider the quotient algebra $L/I_{\lambda}$. For simplicity we
will also denote the restrictions of $R$ and $R^*$ on
$L/I_{\lambda}$ as $R$ and $R^*$ respectively. Since $r$ is a
solution of CYBE on $L$, then by proposition 1 for all $a,b\in L$
(and, consequently, for for all $a,b\in L/I_{\lambda}$):
$$
[R(a),R(b)]-R([a,R(b)])+R([R^*(a),b])=0.
$$

But by the definition of  of the quotient algebra $L/I_{\lambda}$:
$[R^*(a),b]=-[R(a),b]-\lambda[a,b]$  for all $a,b\in L/I_{\lambda}$.
Thus, $R$ is a Rota-Baxter operator of weight $\lambda$ on
$L/I_{\lambda}$.

\begin{corollary}  Let $r=\sum a_i\otimes b_i\in L\otimes L$ be a
solution of CYBE on $L$ such that $r+\tau(r)$ is $L$-invariant. If
$R$ is a Rota-Baxter operator of a weight $\lambda$ then in $L$
there is an ideal $I_{\lambda}$ satisfying $I_{\lambda}\subset
[L,L]$ and $R(I_{\lambda})=0$ and the restriction of the operator
$R$ on the quotient algebra $L/I_{\lambda}$ is a Rota-Baxter
operator of weight $\lambda$ satisfying \eqref{con1}.

Conversely, let $R$ be a Rota-Baxter operator of weight $\lambda$
and for all $a,b\in L$: $R([a,b])+R^*([a,b])=[R(a),b]+[R^*(a),b]$.
Suppose there is an ideal $I_{\lambda}$ in $L$ satisfying
$I_{\lambda}\subset [L,L]$ and $R(I_{\lambda})=0$.  Suppose in
addition that the restriction of the operator $R$ on the quotient
algebra $L/I_{\lambda}$ is a Rota-Baxter operator of weight
$\lambda$ satisfying \eqref{con1}. Then $r$ is a solution of the
CYBE and $r+\tau(r)$ is a nonzero $L$-invariant element.
\end{corollary}

{\bf Proof} Define $I_{\lambda}$ as in \eqref{ideal}. By theorem 1
$R(I_{\lambda})=0$ so $I_{\lambda}$ is $R$-invariant. The rest
follows from theorem 2.

Now suppose that $R$ is a Rota-Baxter operator of weight $\lambda$
and for all $a,b\in L$: $R([a,b])+R^*([a,b])=[R(a),b]+[R^*(a),b]$.
By proposition 1 the element $r+\tau(r)$ is $L$-invariant. Let $x\in
[L,L]$. In the quotient algebra $L/I_{\lambda}$ we have that
$\overline{\theta_{\lambda}(x)}=0$. Thus, $\theta_{\lambda}(x)\in
I_{\lambda}$  and since $R(I_{\lambda})=0$ we have that
$R(\theta_{\lambda})=0$ for all $x\in [L,L]$. By theorem 1 $r$ is a
solution of CYBE on $L$.

We can use the obtained results and specify the main result of
\cite{GME}:

\begin{corollary} Let $L$ be a simple Lie algebra and
$r=\sum\limits_ia_i\otimes b_i\neq 0$ be  a solution of CYBE
\eqref{lieYB} such that $r+\tau(r)\neq 0$ and $r+\tau(r)$ is
$L$-invariant. Then there exists a non-degenerate symmetric
invariant bilinear form $\omega$ on $L$ such that an operator
$R:L\rightarrow L$ defined as
$$
R(a)=\sum\limits_i \omega(b_i,a) a_i
$$
is a Rota---Baxter operator of a nonzero weight $\lambda$ and for
all $a\in L$:
$$
R(a)+R^*(a)+\lambda a=0.
$$
\end{corollary}

{\bf Proof} The existence of the form $\omega$ was proved in
\cite{GME} (theorem 4). Consider the ideal $I_{\lambda}$ defined as
in \eqref{ideal}. Since $L$ is simple, $I_{\lambda}=L$ or
$I_{\lambda}=0$. Since $R$ is a Rota-Baxter operator of weight
$\lambda$, then by theorem 1 $R(I_{\lambda})=0$. If $I_{\lambda}=L$
then $R=0$ that is impassible since $r\neq 0$. Thus, $I_{\lambda}=0$
and by corollary 1 $R(a)+R^*(a)+\lambda a=0$ for all $a\in L$.

\begin{theorem} Let $r=\sum a_i\otimes b_i$ be a non-skew-symmetric
solution of CYBE such that $r+\tau(r)$ is $L$-invariant. Then
operators $R$ and $R^*$ are Rota-Baxter operators of the same
nonzero
 weight $\lambda$ if and only if the commutator ideal $[L,L]$
is the sum of two ideals
$$
[L,L]=I_1\oplus I_2
$$
 such that $R(I_1)=R^*(I_1)=0$ and for all $x\in I_2$:
 $$
 R(x)+R^*(x)+\lambda x=0.
 $$

Conversely, suppose that $R$ and $R^*$ are  Rota-Baxter operators of
the same nonzero weight $\lambda$ such that for all  $a,b\in L$:
$R([a,b])+R^*([a,b])=[R(a),b]+[R^*(a),b]$. Suppose that the
commutator ideal $[L,L]$ is the sum of two ideals $[L,L]=I_1\oplus
I_2$ and the ideals $I_j$ satisfy the same properties as above. Then
$r$ is a solution of the CYBE and $r+\tau(r)$ is a nonzero
$L$-invariant element.
\end{theorem}

{\bf Proof} Suppose $R$ and $R^*$ are Rota-Baxter operators of the
same nonzero weight $\lambda$. By theorem 1,
$R(\theta_{\lambda}([a,b]))=R^*(\theta_{\lambda}[a,b])=0$ for all
$a,b\in L$.

Let $I_1=\{\theta_{\lambda}(a)|\ a\in [L,L]\}$. From proposition 1
it is follows that
$\theta_{\lambda}([a,b])=[\theta_{\lambda}(a),b]=[a,\theta_{\lambda}(b)]$.
Therefore, $I$ is an ideal in $L$. Note, that $R(I_1)=R^*(I_1)=0$.

Let $I_2=\{a\in [L,L]|\ \theta_{\lambda}(a)=0\}$. By proposition 1
$I_2$ is an ideal in L.

If $a\in I_1\cap I_2$, then $R(a)=R^*(a)=0$. But
$0=\theta_{\lambda}(a)=R(a)+R^*(a)+\lambda a=\lambda a$. Thus,
$I_1\cap I_2=0$.

Let $x\in [L,L]$. Then
$\theta_{\lambda}(x-\frac{1}{\lambda}\theta_{\lambda}(x))=\theta_{\lambda}(x)-\theta_{\lambda}(x)=0$.
Thus, $x-\frac{1}{\lambda}\theta_{\lambda}(x)\in I_2$. It means that
$[L,L]=I_1\oplus I_2$.

Conversely, suppose $L=I_1\oplus I_2$ for ideals $I_1$ and $I_2$
satisfying the condition of the theorem. Then
$\theta_{\lambda}([a,b])=0$ if $a\in I_2$ or $b\in I_2$. Let $a,b\in
I_1$. We have
$$
R(\theta_{\lambda}([a,b]))=R(\lambda[a,b])=0=R^*(\theta_{\lambda}([a,b])).
$$
 Thus, $R(\theta_{\lambda}([a,b])=R^*(\theta_{\lambda}([a,b]))=0$ and
 by theorem 1 $R$ and $R^*$ are Rota-Baxter operators of the same nonzero weight
 $\lambda$.

Consider the  inverse statement. For $x\in [L,L]$ we have
$x=i_1+i_2$, where $i_j\in I_j$. It means that
$\theta_{\lambda}(x)=\theta_{\lambda}(i_1)+\theta_{\lambda}(i_2)=\lambda
i_1\in ker(R)$. Therefore, $r$ is a solution of CYBE by theorem 1
and by proposition 1 $r+\tau(r)$ is $L$-invariant.

\begin{corollary}
Let $(L,\omega)$ be a quadratic Lie algebra and $R$ be a Rota-Baxter
operator of a nonzero weight $\lambda$. Suppose that for all $a\in
L$
$$
R(a)+R^*(a)+\lambda a=0.
$$
Then the corresponding tensor $r$ is a solution of CYBE and
$r+\tau(r)$ is $L$-invariant.
\end{corollary}

\begin{example}
Consider $(L_2,\omega)$ --- quadratic lie algebra from example 3.
Direct computations show:
$$
R(a)+R^*(a)=a+a\overline t, \ \ R(a\overline t)+R^*(a\overline
t)=a\overline t,
$$
for all $a\in L$. Thus, for all $a,b\in L$:
$$
[R(a)+R^*(a),b]=[a+a\overline t,b]=R([a,b])+R^*([a,b]),
$$
$$
[R(a)+R^*(a),b\overline t]=[a,b]\overline t=R([a,b\overline
t])+R^*([a,b\overline t]),
$$
$$
[R(a\overline t)+R^*(a\overline t),b\overline t])=0=R([a\overline
t,b\overline t])+R^*([a\overline t,b\overline t]).
$$
This means that for the corresponding tensor $r=\frac{1}{4}h\otimes
h\overline{t}+x\otimes y\overline{t}+y\otimes x\overline{t}$ the
element $r+\tau(r)$ is $L$-invariant.
\end{example}

At the end of the section, we will specify the results in the case
of algebraically closed field.

\begin{theorem} Let $(L,\omega)$ be a quadratic Lie algebra over an algebraically closed field $F$ and let $r=\sum a_i\otimes b_i\in
L\otimes L$ be a solution of CYBE equation such that $r+\tau(r)$ is
$L$-invariant. Then operators $R$ and $R^*$ are  Rota-Baxter
operator of the same nonzero weight $\lambda$ if and only if $L$ is
the sum of ideals $L=I_1\oplus I_2$ such that for all $a\in I_1$:
$$
R(a)+R^*(a)+\lambda a\in Z(L)
$$
and for all $x\in I_2$ and $y\in L$:
$$
R([x,y])=R^*([x,y])=0.
$$

Conversely, suppose that  $R$ is a Rota-Baxter operator of a nonzero
weight $\lambda$, such that for all for all  $a,b\in L$:
$R([a,b])+R^*([a,b])=[R(a),b]+[R^*(a),b]$. If $L$ is a direct sum of
two ideals $L=I_1\oplus I_2$  and the ideals $I_j$ satisfy the same
properties as above, then $r$ is a solution of CYBE equation and
$r+\tau(r)$ is $L$ invariant.
\end{theorem}

{\bf Proof} Let $R$ be a Rota-Baxter operator of weight $\lambda$.
Consider the linear map $\theta=\theta_0=R+R^*:L\mapsto L$. Since
the ground field is algebraically closed, $L$ is the sum of root
subspaces
 $$
 L=L(\alpha_1)\oplus\ldots\oplus L(\alpha_k),
 $$
 where $\{\alpha_1,\ldots,\alpha_k\}$ is the spectrum of $\theta$ and $L(\alpha_i)=\{v\in L|\ (\theta-\alpha_i  id)^{n_i}(v)=0\ \text{for some $n_i$}\}$ -- is the root subspace corresponding to the eigenvalue $\alpha_i$. Since $\theta_{\lambda}([a,b])=[\theta_{\lambda}(a),b]$, $L(\alpha_i)$ is an ideal in $L$.

 Put $I_1=L(-\lambda)$ if $-\lambda$ is an eigenvalue of $\theta$ and $I_1=0$ otherwise. Then $L=I_1\oplus I_2$, where $I_2$ is efore, every root subspace is an ideal in $L$ and, consequently, $I_1$ and $I_2$ are ideals in $L$.

Take $a\in I_1$. By the definition of $I_1$:
$\theta_{\lambda}^k(a)=0$ for some $k\in \mathbb N$. First consider
the case when $a=[x,y]$ for some $x,y\in I_1$. We have
$$
\theta_{\lambda}^2([x,y])=\theta_{\lambda}(R([x,y])+R^*([x,y])+\lambda[x,y])=
$$
$$
=(R+R^*)([R(x),y]+[R^*(x),y]+\lambda[x,y])+\lambda\theta_{\lambda}([x,y])=\lambda\theta_{\lambda}([x,y]).
$$

Therefore,
$$0=\theta_{\lambda}^k([x,y])=\lambda^{k-1}\theta_{\lambda}([x,y]).$$
And since $\lambda\neq 0$ we finally obtain that
$\theta_{\lambda}([x,y])=R([x,y])+R^*([x,y])+\lambda[x,y]=0$ for all
$x,y\in I_1$.

Now let $a$ be an arbitrary element in $I_1$. It is obvious that
$[\theta_{\lambda}(a),I_2]=0$. If $b\in I_1$, then
$[\theta_{\lambda}(a),b]=\theta_{\lambda}([a,b])=0$. Thus,
$\theta_{\lambda}(a)\in Z(L)$ for all $a\in I_2$.

Consider the ideal $I_2$. By the properties of root subspaces, $I_2$
is $\theta_{\lambda}$-invariant and the restriction of
$\theta_{\lambda}$ on $I_2$ is invertible. Thus, for every $x\in
I_2$: $x=\theta_{\lambda}(x')$ for some $x'\in I_2$.

Let $x,y\in I_2$. Then
$$
R([x,y])=R([\theta_{\lambda}(x'),y])=R(\theta_{\lambda}([x',y]))=0
$$
by theorem 1. Similarly, $R^*([x,y])=0$.

The rest statements of the theorem follows from theorem 3. Indeed,
if $L=I_1\oplus I_2$ , then $[L,L]=[I_1,I_1]\oplus[I_2,I_2]$. Since
$\theta_{\lambda}(a)\in Z(L)$ for all $a\in I_1$, then
$\theta_{\lambda}([a,b])=[\theta_{\lambda}(a),b]=0$ for all $a,b\in
I_1$. Also, it is obvious that for all $x\in [I_2,I_2]$
$R(x)=R^*(x)=0$ and we may use theorem 3 to prove the rest two
statements of the theorem.




\begin{thebibliography}{1}

\bibitem{Br} Baxter G., \emph{An analytic problem whose solution follows from a simple algebraic identity}, Pacific J. Math., 10 (1960),
731--742. 

\bibitem{Atk} Atkinson, F.V., \emph{Some aspects of BaxterТs functional equation}, J. Math. Anal. Appl.,
7  (1963), 1Ц30. 

\bibitem{Rota}  Rota G.C., \emph{Baxter algebras and combinatorial identities I and II}, Bull. Amer. Math.
Soc., 75 (1969),  325--334. 

\bibitem{Miller} Miller J.B., \emph{Some properties of Baxter operators}, Acta Math. Acad. Sci.
Hungar.,  17 (1966), 387--400. 

\bibitem{Car} Cartier P., \emph{On the structure of free Baxter algebras}, Adv. Math., 9 (1972),
253--265. 

\bibitem{Guo} Guo L., \emph{An Introduction to Rota---Baxter Algebra}, Surveys of Modern
Mathematics,  Somerville, MA: International Press; Beijing: Higher
education press,   4 (2012). 


\bibitem{BD} Belavin A.A., Drinfeld V.G., \emph{Solutions of the classical Yang---Baxter equation for simple Lie algebras}, Funct. Anal.
Appl., 16:3 (1982),  159Ц180. 

\bibitem{STS} Semenov-Tyan-Shanskii M.A., \emph{What a classical r-matrix is},  Funct. Anal.
Appl., 17:4 (1983),  259--272. 

\bibitem{GME} Goncharov M. E. On Rota-Baxter operators of non-zero weight arisen from the solutions of the classical Yang-Baxter equation,
 Sib. El. Math. Rep., 14 (2017) 1533-1544

 \bibitem {Drinf}
Drinfeld V.G., \emph{Hamiltonian structures on Lie groups, Lie
bialgebras and the geometric meaning of the classical Yang---Baxter
equation}, Sov, Math, Dokl, 27 (1983), 68--71.

\bibitem {ANQ}
Anquela J.A., Cortes T., Montaner F., \emph{Nonassociative
Coalgebras}, Comm. Algebra., 22:12 (1994),  4693--4716.

\bibitem {Zhelyabin98}
Zhelyabin V.N., \emph{Jordan bialgebras of symmetric elements and
Lie bialgebras}, Sib. Math. J.,  39:2 (1998), 261--276. 


\bibitem {Zhelyabin}
Zhelyabin V.N., \emph{On a class of Jourdan D-bialgebras}, St.
Petersburg Mathematical Journal,  11:4 (2000), 589--609.

\bibitem{versh} Vershinin V.V., \emph{On Poisson--Malcev Structures}, Acta Applicandae
Mathematicae, 75 (2003),  281--292.


\bibitem {Aquiar}
Aguiar M., \emph{On the associative analog of Lie bialgebras},
Journal of Algebra, 244 (2001), 492--532.

\bibitem{Polishchuk}  Polishchuk A., \emph{Clasic Yang---Baxter Equation and the
A-constraint}, Advances in Mathematics,  168:1 (2002), 56--96.


\bibitem{Gme} Goncharov M.E., \emph{The classical Yang---Baxter equation on alternative algebras:
The alternative D-bialgebra structure on Cayley--Dickson matrix
algebras},  Sib. Math. J.,  48:5 (2007), 809--823.

\bibitem {GMM} Goncharov M.E., \emph{Structures of Malcev Bialgebras on a Simple Non-Lie Malcev Algebra},  Commun. Algebra., 40:8 (2012),
3071--3094.

\bibitem{BGN} Bai C., Guo L. and  Ni X., \emph{Nonabelian generalized Lax pairs, the classical Yang-Baxter equation and Post Lie
algebras}, Comm. Math. Phys., 297, (2010), 553-596.

\bibitem{BGN1} Bai C., Guo L. and  Ni X., \emph{Generalizations of the classical Yang-Baxter equation and
${O}$-operators}, Journal of Math. Phys., 52 (2011), 063515.

\bibitem{BajoBen} Bajo I.,  Benayadi s. \emph{Lie algebras admitting a unique quadratic
structure}, Commun. Algebra. 25:9 (1997), 2795--2805.

\end{thebibliography}
\end{document}